\documentclass[11pt,a4paper]{article}
\usepackage[latin1]{inputenc}
\usepackage[english]{babel}
\usepackage[T1]{fontenc}

\usepackage{amsmath}
\usepackage{amsthm}
\usepackage{amsfonts}
\usepackage{amssymb}
\usepackage{amscd}
\usepackage{eucal}
\usepackage{upref}
\usepackage{fancyhdr}
\usepackage[all]{xy}

\renewcommand{\mathbf}{\boldsymbol}

\theoremstyle{plain}
\newtheorem{theorem}{Theorem}[section]
\newtheorem{lemma}[theorem]{Lemma}
\newtheorem{proposition}[theorem]{Proposition}

\numberwithin{equation}{section}

\theoremstyle{definition}

\newtheorem{example}[theorem]{Example}

\theoremstyle{remark}

\DeclareMathOperator{\MI}{I}
\DeclareMathOperator{\End}{End}
\DeclareMathOperator{\ad}{ad}
\DeclareMathOperator{\Ric}{Ric}

\DeclareMathOperator{\tr}{tr}
\DeclareMathOperator{\image}{Im}

\DeclareMathOperator{\Hess}{Hess}
\DeclareMathOperator{\Ad}{Ad}
\DeclareMathOperator{\Hom}{Hom}
\DeclareMathOperator{\Sym}{Sym}

\let\GAMMA=\Gamma
\renewcommand{\Gamma}{\mathit{\GAMMA}}
\let\DELTA=\Delta
\renewcommand\Delta{\mathit{\DELTA}}
\let\THETA=\Theta
\renewcommand{\Theta}{\mathit{\THETA}}
\let\LAMBDA=\Lambda
\renewcommand{\Lambda}{\mathit{\LAMBDA}}
\let\XI=\Xi
\renewcommand{\Xi}{\mathit{\XI}}
\let\PI=\Pi
\renewcommand{\Pi}{\mathit{\PI}}
\let\SIGMA=\Sigma
\renewcommand{\Sigma}{\mathit{\SIGMA}}
\let\PHI=\Phi
\renewcommand{\Phi}{\mathit{\PHI}}
\let\PSI=\Psi
\renewcommand{\Psi}{\mathit{\PSI}}
\let\OMEGA=\Omega
\renewcommand{\Omega}{\mathit{\OMEGA}}

\title{Initial value problem for cohomogeneity one gradient Ricci solitons}

\author{Maria Buzano\\
{\small Wolfson College, Oxford OX2 6UD, UK}\\
{\small \texttt{buzano@maths.ox.ac.uk}}}

\begin{document}

\maketitle

\begin{abstract}
Consider a smooth manifold $M$. Let $G$ be a compact Lie group which acts on $M$ with cohomogeneity one. Let $Q$ be a singular orbit for this action. We study the gradient Ricci soliton equation $\Hess(u)+\Ric(g)+\frac{\epsilon}{2}g=0$ around $Q$. We show that there always exists a solution on a tubular neighbourhood of $Q$ for any prescribed $G$-invariant metric $g_Q$ and shape operator $L_Q$, provided that the following technical assumption is satisfied: if $P=G/K$ is the principal orbit for this action, the $K$-representations on the normal and tangent spaces to $Q$ have no common sub-representations. We also show that the initial data are not enough to ensure uniqueness of the solution, providing examples to explain this indeterminacy. This work generalises the papaer "The initial value problem for cohomogeneity one Einstein metrics" of 2000 by J.-H. Eschenburg and McKenzie Y. Wang to the gradient Ricci solitons case.
\end{abstract}


\small \textbf{MSC}: Differential geometry

\section{Introduction}
A \emph{Ricci soliton} is a triple $(M^n,\widehat{g},\widehat{X})$, where $(M,\widehat{g})$ is a complete Riemannian manifold of dimension $n$ and $\widehat{X}$ is a vector field on $M$, which satisfies the equation
\begin{equation}\label{eqn:0.2}
\Ric(\widehat{g})+\frac{1}{2}\mathcal L_{\widehat{X}}\widehat{g}+\frac{\epsilon}{2}\widehat{g}=0,
\end{equation}
where $\epsilon$ is some real constant and $\Ric(\widehat{g})$ is the Ricci tensor of the metric $\widehat{g}$. 

Equation \eqref{eqn:0.2} may be written also as
\begin{equation}\label{eqn:0.3}
\Ric(\widehat{g})+\delta^*\widehat{\omega}+\frac{\epsilon}{2}\widehat{g}=0,
\end{equation}
where $\widehat{\omega}=\widehat{X}^\flat$ and $\delta^*$ is the symmetrized covariant derivative. If the dual form of $X$ is exact, i.e. if there exists a smooth function $u$ on $M$ such that $\widehat{X}^\flat=du$, the Ricci soliton is called a \emph{gradient Ricci soliton} and equation \eqref{eqn:0.3} takes the following form
\begin{equation}\label{eqn:0.4}
\Ric(\widehat{g})+\Hess(u)+\frac{\epsilon}{2}\widehat{g}=0,
\end{equation}
where $\Hess(u)$ is the Hessian of $u$.

Ricci solitons are important for many reasons. Firstly, they generate particular solutions to the Ricci flow equation. In fact, considering the one-parameter family of vector fields on $M$ given by
\begin{equation*}
X(t)=\frac{\widehat{X}}{1+\epsilon t},
\end{equation*}
and integrating it to a one-parameter family $\varphi(t)$ of diffeomorphisms of $M$, we have that the one-parameter family of Riemannian metrics
\begin{equation*}
g(t)=(1+\epsilon t)\varphi(t)^*\widehat{g}
\end{equation*}
evolves under the Ricci flow equation
\begin{equation*}
\begin{cases}
&\frac{\partial g(t)}{\partial t}=-2\Ric(g(t)),\\
&g(0)=\widehat{g},
\end{cases}
\end{equation*}
for $t\in[0,T)$, where $T\in(0,+\infty)$. Secondly, the notion of gradient Ricci soliton has motivated the discovery of monotonicity formulas for the Ricci flow (see \cite{Perelman:1}), which have many geometric applications. Thirdly, Ricci solitons often appear as limits of dilations of singularities in the Ricci flow (see \cite{Hamilton:1}). Finally, they are generalisations of Einstein metrics. In fact, if we take $\widehat{X}$ to be the zero vector field in \eqref{eqn:0.2} or if we take $u$ to be constant in \eqref{eqn:0.4}, we obtain the Einstein condition for the metric $\widehat{g}$, which is
\begin{equation*}
\Ric(\widehat{g})+\frac{\epsilon}{2}\widehat{g}=0.
\end{equation*}
Therefore, it is natural to ask whether techniques used to produce Einstein metrics can also produce Ricci solitons. For example, we may consider spaces with a certain amount of symmetry. By \cite{Petersen-Wylie:1}, we have that the maximal amount of symmetry on a nontrivial (i.e. non Einstein) gradient Ricci soliton is given by a cohomogeneity one action. This means that the symmetry group of the manifold acts on it in such a way that the orbits with maximal dimension, which are called \emph{principal orbits}, have codimension one. In this context, orbits with lower dimension are called \emph{singular orbits} or \emph{special orbits}. With this kind of action, the orbit space is one-dimensional and can be an interval (open, closed or semi-open), $\mathbb R$ itself or the circle $S^1$, depending on the number of singular orbits (cf. \cite{Bergery:1}).

Many examples of cohomogeneity one gradient solitons have been constructed; e.g., the Bryant soliton (see \cite{Ivey:1}, \cite{Chow:1}), the cigar soliton (\cite{Hamilton:3}), the $U(n)$-symmetric soliton on $\mathbb C^n$ discovered by Cao (\cite{Cao:2}, \cite{Cao:3}) and the K\"{a}hler generalisations \cite{Koiso:1}, \cite{Yang:1}, \cite{Dancer-Wang:2}, \cite{Eminenti-Mantegazza-LaNave:1} and \cite{Feldman-Ilmanen-Knopf:1}.

In \cite{Wang-Eschenburg:1}, the authors consider local existence and uniqueness of a smooth $G$-invariant Einstein metric around a singular orbit in the cohomogeneity one setting. They prove, under a technical assumption, that, given any $G$-invariant Riemannian metric and any shape operator in a neighbourhood of a singular orbit $Q$, there always exists an invariant Einstein metric around $Q$ with any prescribed sign of the Einstein constant. We generalise the Einstein case to the case of gradient Ricci solitons. In particular, we study the existence and uniqueness of $G$-invariant gradient Ricci solitons around the  singular orbit $Q$. Our main theorem is the following.
\begin{theorem}
Let $G$ be a compact Lie group acting on a connected Riemannian manifold $(M,\widehat{g})$ with cohomogeneity one and by isometries of $\widehat{g}$. Let $Q=G/H$ be a singular orbit of codimension $k+1$, $k\geq1$. So, $H$ acts linearly with cohomogeneity one on $V=\mathbb R^{k+1}$ and the Lie algebra of $G$ splits as $\mathfrak g=\mathfrak h\oplus\mathfrak p_-$, where $\mathfrak h$ is the Lie algebra of $H$. Let $v_0\in S^k$ have isotropy group $K$ with respect to the $H$-action. Then, $G/K$ is a principal orbit for the action.

Assume that $V$ and $\mathfrak p_-$ have no irreducible common factors as $K$- representations. Then, given any $G$-invariant metric $g_Q$ on $Q$ and any shape operator $L_1:NQ\to\Sym^2(T^*Q)$, where $NQ=G\times_H V$ is the normal bundle over $Q$, there exists a $G$-invariant gradient Ricci soliton on some open disk bundle of $NQ$.
\end{theorem}
We write the Ricci soliton condition around $Q$ as an initial value problem with initial data given by a $G$-invariant metric and shape operator $L_1$ around $Q$. We note that the smoothness condition of the metric implies that $\tr{(L_1)}=0$, so $Q$ must be a minimal submanifold in $M$. As we are working around a singular orbit, we only need to solve the gradient Ricci soliton equation in the directions tangent and orthogonal to the orbits. Moreover, we can write the initial value problem as a system of ordinary non-linear differential equations with a singular point at the origin. We solve this system using the same technique as in \cite{Wang-Eschenburg:1}, which consists of applying the method of asymptotic series to find a solution and then showing that this solution is in fact a smooth $G$-invariant gradient Ricci soliton. We can always show existence, but the initial data given are not sufficient to ensure uniqueness and the indeterminacy of the problem, which is the same as in the Einstein case, is always finite and can be computed using representation theory. Finally, the technical assumption about the irreducible summands of $V$ and $\mathfrak p_-$ is motivated by \cite{Wang-Eschenburg:1} and, as the authors explain (cf. Remark 2.7 of \cite{Wang-Eschenburg:1}), it appears quite natural in the context of the Kaluza-Klein construction.

\section{The cohomogeneity one Ricci soliton equation}
In this section, following the approach and notation of \cite{Wang-Eschenburg:1} and \cite{Dancer-Wang:2}, we recall the Ricci soliton equation in the cohomogeneity one setting. 

Let $(M,\widehat{g})$ be a connected Riemannian manifold of dimension $n+1$ and let $G$ be a compact Lie group which acts on $M$ by isometries and with cohomogeneity one.

Now, let us choose a unit speed geodesic
\begin{equation*}
\gamma:I\longrightarrow M,
\end{equation*}
where $I\subset\mathbb R$ is an open interval, and take $\gamma$ such that it intersects all the principal orbits orthogonally. Then, it is possible to define an equivariant diffeomorphism
\begin{equation*}
\begin{aligned}
\Phi:\,&I\times G/K\longrightarrow M_0\subset M,\\
&(t,g\cdot K)\longmapsto g\cdot\gamma(t),
\end{aligned}
\end{equation*}
where $K$ is the isotropy group of $\gamma(t)$ with respect to the $G$-action. Hence, $\Phi(t,G/K)$ is the principal orbit $P_t$ passing through $\gamma(t)$ and $M_0$ is an open dense subset in $M$ which is the union of all the principal orbits.

Every orbit $P_t$ is naturally equipped with a $G$-invariant Riemannian metric, which depends on $t$. Hence, through $\Phi$ we obtain a family $g(t)$ of $G$-invariant metrics on the homogeneous space $P$, where $P$ denotes an abstract copy of the principal orbit $G/K$. Moreover, the map $\Phi$ sends $I\times\{p\}$, with $p\in P$, to the geodesic through $p$ orthogonal to the principal orbits, and we have that the canonical parametrisation of $I$ corresponds to the arclength parametrisation of the geodesic. Then, if we pullback the metric $\widehat{g}$ on $M$ through $\Phi$ we get
\begin{equation*}
\Phi^*(\widehat{g})=d t^2+g_t,
\end{equation*}
where $g_t$ is a family of $G$-invariant Riemannian metrics on $P$. 

Let $\widehat{\nabla}$ and $\widehat{\Ric}$ denote the Levi-Civita connection and the Ricci tensor of the manifold $(M,\widehat{g})$, respectively. Let $\nabla_t$ and $\Ric_t$ denote the Levi-Civita connection and the Ricci tensor of $(P_t,g_t)$, respectively. Let $L_t$ be the shape operator on $P_t$ defined by
\begin{equation*}
L_t(X)=\widehat{\nabla}_XN,
\end{equation*}
where $X$ is a vector field on $P_t$ and $N=\Phi_*(\frac{\partial}{\partial t})$ is a unit normal $G$-invariant vector field along $P_t$ such that $\widehat{\nabla}_NN=0$. We then have a family $(L_t)_{t\in I}$ of $G$-invariant, $g_t$-symmetric endomorphisms of the tangent space of $P$. In particular, the trace of $L_t$ is constant along $P_t$. Moreover, the following equality holds
\begin{equation}\label{eqn:2.1}
\dot{g}_t(X,Y)=2g_t(L_t(X),Y),
\end{equation}
for every pair of vector fields $X,Y$ on $P_t$ and for every $t\in I$. 

Consider now the Ricci soliton equation for $(M,\widehat{g},\widehat{\omega})$
\begin{equation}\label{eqn:2.6}
\widehat{\Ric}(\widehat{g})+\widehat{\delta}^*\widehat{\omega}+\frac{\epsilon}{2}\widehat{g}=0. 
\end{equation}
First of all, we can take the form $\widehat{\omega}$ to be $G$-invariant, or, if we are dealing with gradient Ricci solitons, we can take the potential function to be $G$-invariant (cf. \cite{Dancer-Wang:2}, p. 4). Hence, if we consider the pull-back of $\widehat{\omega}$ through $\Phi$, we obtain
\begin{equation}\label{omega}
\Phi^*\widehat{\omega}=\xi(t)d t+\omega_t,
\end{equation}
where $\xi=\xi(t)$ is a function on $I$ and $\omega_t$ is a one-parameter family of $G$-invariant one-forms on $P$. 

Let $\mathcal X(P_t)$ denote all the vector fields on $P_t$. We then have the following proposition.
\begin{proposition}[\cite{Dancer-Wang:2}]
Let $(M^{n+1},\widehat{g})$ be a connected Riemannian manifold which admits a cohomogeneity one action by a compact Lie group $G$ of isometries of $\widehat{g}$. Let $\widehat{\omega}$ be a $G$-invariant one-form on $M$. Under the parametrisation induced by a unit speed geodesic orthogonal to the principal orbits, the Ricci soliton equation for $\widehat{g}$ and the vector field dual to $\widehat{\omega}$ is given by
\begin{align*}
&-(\delta^{\nabla_t}L_t)^\flat-d(\tr(L_t))+\frac{1}{2}\dot{\omega}_t-\omega_t\circ L_t=0,\\
&-\tr(\dot{L}_t)-\tr(L^2_t)+\dot{\xi}(t)+\frac{\epsilon}{2}=0,\\ 
\begin{split}
&\Ric_t(X,Y)-\tr(L_t)g_t(L_t(X),Y)-g_t(\dot{L}_t(X),Y)\\ 
&+\xi(t)g_t(L_t(X),Y)+\delta^*_t\omega_t(X,Y)+\frac{\epsilon}{2}g_t(X,Y)=0,
\end{split}
\end{align*}
for all $X,Y\in\mathcal X(P_t)$ and $t\in I$, where, viewing $L_t$ as an endomorphism of $\mathcal X(P_t)$, the operator $\delta^{\nabla_t}:\mathcal X^*(P_t)\otimes\mathcal X(P_t)\to\mathcal X(P_t)$ is the codifferential.

Conversely, if $g_t$ and $\omega_t$ are one-parameter families of metrics and one-forms on $P_t$, respectively, and $\xi=\xi(t)$ is a smooth function on $I$ such that the above system is satisfied with $L_t$ defined by $\dot{g}_t(X,Y)=2g_t(L_t(X),Y)$ for all $X,Y\in\mathcal X(P_t)$, then $\widehat{g}=d t^2+g_t$ and $\widehat{\omega}=\xi(t)d t+\omega_t$ give a local Ricci soliton on $M_0$.
\end{proposition}
If we are looking for gradient Ricci solitons, that is when there exists a $G$-invariant smooth function $u$ such that $\widehat{\omega}=d u$, equation \eqref{omega} becomes
\begin{equation*}
\Phi^*\widehat{\omega}=\dot{u}(t)d t,
\end{equation*}
and the Ricci soliton equation in the cohomogeneity one setting is equivalent to the following system
\begin{align}
\label{eqn:2.11} &-(\delta^{\nabla_t}L_t)^\flat-d(\tr(L_t))=0,\\
\label{eqn:2.12} &-\tr(\dot{L}_t)-\tr(L^2_t)+\ddot{u}(t)+\frac{\epsilon}{2}=0,\\ 
\label{eqn:2.13}\begin{split}
&\Ric_t(X,Y)-\tr(L_t)g_t(L_t(X),Y)-g_t(\dot{L}_t(X),Y)+\dot{u}(t)g_t(L_t(X),Y)\\
&+\frac{\epsilon}{2}g_t(X,Y)=0,\end{split}
\end{align}
for all $X,Y\in \mathcal X(P_t)$ and $t\in I$, where $u(t)(p)=u\circ\Phi(t,p)$, for all $p\in P_t$.

\section{Smoothness of tensors around a singular orbit}\label{section 3}
In this section, following \cite{Wang-Eschenburg:1}, we will discuss briefly the smoothness criterion for the metric $\widehat{g}$ and the one-form $\widehat{\omega}$,  in the case when there is a special orbit.

Let $(M,\widehat{g})$ be a connected $(n+1)$-dimensional Riemannian manifold and $G$ a compact Lie group which acts on $M$ by isometries of $\widehat{g}$ and with cohomogeneity one. Let $Q=G\cdot q$ be a singular orbit of codimension $k+1$, with $k\geq1$, with isotropy group $H=G_q$. As $Q=G/H$ is a homogeneous space, the Lie algebra $\mathfrak g$ of $G$ decomposes in the following way:
\begin{equation}\label{eqn:3.1}
\mathfrak g=\mathfrak h\oplus\mathfrak p_-,
\end{equation}
where $\mathfrak h$ is the Lie algebra of $H$ and $\mathfrak p_-$ is the $\Ad(H)$-invariant complement of $\mathfrak h$ in $\mathfrak g$, which can be identified with the tangent space $T_qQ$.

Let $V=T_qM/T_qQ\backsimeq\mathbb R^{k+1}$ be the normal space at $q$ of $Q$, on which $H$ acts linearly with cohomogeneity one, i.e. it acts transitively on the sphere $S^k=H/K$, where $K\subset H$. We can identify a tubular neighbourhood of $Q$ with the total space of the normal bundle $NQ=G\times_HV$ of $Q$. We have that 
\begin{equation*}
T(NQ)_{|_V}=V\times(V\oplus\mathfrak p_-).
\end{equation*}
In fact, the tangent space to $NQ$ splits into orthogonal and vertical parts:
\begin{equation*}
T(NQ)=\pi^*NQ\oplus\pi^*TQ,
\end{equation*}
where $\pi:NQ\rightarrow Q$ is the bundle map, and the two pull-back bundles are trivial on $V$, which can be viewed as the fibre of $NQ$ over $q\in Q$. Hence, a smooth $G$-invariant symmetric bilinear form $a$ is determined by an $H$-equivariant smooth map
\begin{equation*}
a:V\longrightarrow \Sym^2(V\oplus\mathfrak p_-).
\end{equation*} 

Let $W$ be the vector space of all smooth $H$-equivariant maps $L:S^k\to\Sym^2(V\oplus\mathfrak p_-)$. We then have that, for $v_0\in S^k$, the evaluation map 
\begin{equation*}
\begin{aligned}
\text{ev}:\,&W\longrightarrow \Sym^2(V\oplus\mathfrak p_-)^K,\\
&L\longmapsto\text{ev}(L)=L(v_0),
\end{aligned}
\end{equation*}
is a linear isomorphism. Here, $\Sym^2(V\oplus\mathfrak p_-)^K$ denotes the elements of $\Sym^2(V\oplus\mathfrak p_-)$ which are $K$-invariant, where $K=H_{v_0}$. Let $W_m$ be the subspace of $W$ consisting of all maps which are restrictions to $S^k$ of $H$-equivariant homogeneous polynomials of degree $m$. We then have a necessary and sufficient condition for $a$ to be smooth.
\begin{lemma}[\cite{Wang-Eschenburg:1}]\label{Lemma:3.1}
Let $t\mapsto a_t$, where $a_t: S^k\to\Sym^2(V\oplus\mathfrak p_-)^K$ for all $t\in[0,\infty)$, be a smooth curve, i.e. at zero the right-hand derivatives of all orders exist and are continuous from the right. Let $\sum_p{a_pt^p}$ be its Taylor expansion at zero. Then the map $a$ defined by 
\begin{equation*}
\begin{aligned}
a:V&\backslash\{0\}\longrightarrow\Sym^2(V\oplus\mathfrak p_-),\\
&v\longmapsto a(v)=a_{|v|}\left(\frac{v}{|v|}\right)
\end{aligned}
\end{equation*}
can be extended smoothly at zero if and only if $a_p\in\text{ev}(W_p)$ for all $p\geq0$.
\end{lemma}
Motivated by \cite{Wang-Eschenburg:1}, we now assume that the representations of $K$ on $\mathfrak p_-$ and $V$ have no irreducible common factors. As a consequence of this, we have that 
\begin{equation}\label{assumption}
\Sym^2(V\oplus\mathfrak p_-)^K=\Sym^2(V)^K\oplus \Sym^2(\mathfrak p_-)^K,
\end{equation}
and each $W_m$ splits as $W_m^+\oplus W_m^-$, where the polynomials in $W_m^+$ take values in $\Sym^2(V)$ and the ones in $W_m^-$ take values in $\Sym^2(\mathfrak p_-)$.

The smoothness criterion for $\widehat{\omega}$ is obtained essentially in the same way as for the metric $\widehat{g}$. Under the above assumption, we have that on a tubular neighbourhood around $Q$, $\widehat{\omega}$ is determined by an $H$-equivariant map 
\begin{equation*}
\widehat{\omega}:V\longrightarrow V^*\oplus\mathfrak p_-^*.
\end{equation*}
In this case, $W_m$ is defined as the space of $H$-equivariant maps $L:V\to V^*\oplus\mathfrak p_-^*$ which are restrictions to the unit sphere $S^k$ of homogeneous polynomials of degree $m$. The necessary and sufficient condition for $\widehat{\omega}$ of the form \eqref{omega} to be smooth is that its pth Taylor coefficient, viewing $\widehat{\omega}$ as function of $t$, lives in $\text{ev}(W_p)$, for all $p\geq0$.

Finally, if we consider gradient Ricci solitons with potential function $u$, in the case of a special orbit the smoothness criterion for $\widehat{\omega}=d u$ implies that $u(t)$ must be even in $t$. In fact, around zero, $u$ is given by
\begin{equation}\label{eqn:3.3}
u(t)=\sum_{t=0}^\infty{\frac{u_p}{p!}t^p},
\end{equation}
where
\begin{equation*}
u_p=\frac{d^p}{dt^p}u(t)\Big\vert_{t=0}.
\end{equation*}
The smoothness condition implies that $u_p$ must be a homogeneous polynomial of degree $p$ on the sphere $S^k$, on which $H$ acts transitively. Moreover, we can take $u$ to be $H$-invariant. We now show that $u_p=0$ if $p$ is odd. Given $x\in S^k$, there exists $h\in H$ such that $h\cdot x=-x\in S^k$. If $p$ is odd we have that
\begin{equation*}
-u_p(x)=u_p(-x)=u_p(h\cdot x)=u_p(x)\Longrightarrow u_p(x)=0.
\end{equation*}
Hence, if $p$ is odd, $u_p(x)=0$, for all $x\in S^k$. This implies that $u(t)$ given by \eqref{eqn:3.3} is even in t.

\section{Initial value problem for gradient Ricci solitons around a singular orbit}\label{sec:4}
First of all note that by Proposition 2.17 of \cite{Dancer-Wang:2}, if we are looking for gradient Ricci solitons in the case when there is a singular orbit, instead of considering the system \eqref{eqn:2.11}-\eqref{eqn:2.13}, we can consider the following system
\begin{align}
\label{eqn:4.1} &\frac{d^3}{d t^3}u(t)+\tr(L_t)\ddot{u}(t)+\tr(\dot{L}_t)\dot{u}(t)-2\ddot{u}(t)\dot{u}(t)-\epsilon\dot{u}(t)=0\\
\label{eqn:4.2} \begin{split}
&\Ric_t(X,Y)-\tr(L_t)g_t(L_t(X),Y)-g_t(\dot{L}_t(X),Y)+\dot{u}(t)g_t(L_t(X),Y)\\
&+\frac{\epsilon}{2}g_t(X,Y)=0,
\end{split}
\end{align}
together with \eqref{eqn:2.1}. Note that equation \eqref{eqn:4.1}, which can be viewed as an equation in $\dot{u}(t)$, is the first integral which arises from the contracted second Bianchi identity and was observed in \cite{Ivey:1} (p. 242) and more generally in \cite{Cao:1} (p. 123) and \cite{Hamilton:1} (pp. 84-85). Moreover, the smoothness condition on the function $u$ implies that
\begin{equation*}
\dot{u}(0)=0.
\end{equation*}
Using the Ricci endomorphisms $r_t$ on $P_t$, defined by
\begin{equation*}
\Ric_t(X,Y)=g_t(r_t(X),Y),
\end{equation*}
for all $X,Y\in\mathcal X(P_t)$, equation \eqref{eqn:4.2} becomes
\begin{equation}\label{eqn:4.3}
r_t-\tr(L_t)L_t-\dot{L}_t+\dot{u}(t)L_t+\frac{\epsilon}{2}\MI=0,
\end{equation}
where $\MI$ is the identity matrix. 

Let $\mathfrak h$ in \eqref{eqn:3.1} decompose in the following way
\begin{equation*}
\mathfrak h=\mathfrak k\oplus\mathfrak p_+,
\end{equation*}
and let $\mathfrak p=\mathfrak p_+\oplus\mathfrak p_-$, so that
\begin{equation*}
\mathfrak g=\mathfrak k\oplus\mathfrak p.
\end{equation*}
Hence, $\mathfrak p$ is the tangent space at a point to the principal orbit $G/K$, while $\mathfrak p_+$ and $\mathfrak p_-$ can be identified with the tangent spaces to $H/K$ and $G/H=Q$, respectively.

As we are working around the singular orbit $Q$, by assumption \eqref{assumption}, we have that $g_t$ and $L_t$ split in $+$ and $-$ parts. Hence, following \cite{Wang-Eschenburg:1}, let us choose $x(t),\eta(t)\in\End(\mathfrak p)^K$ preserving the splitting of $\mathfrak p$ and such that 
\begin{align*}
&g_t=t^2x_+(t)\oplus x_-(t),\\
&L_t=\left(\frac{1}{t}\MI_++\eta_+(t)\right)\oplus\eta_-(t),
\end{align*}
with initial conditions given by
\begin{align*}
&x(0)=\MI,\\
&\eta_+(0)=0\quad\text{and}\quad\eta_-(0)=L_1(v_0),
\end{align*}
where $L_1$ is the shape operator of the singular orbit $Q$, which is an $H$-equivariant linear map from $V\longrightarrow\Sym^2(\mathfrak p_-)$. We will now drop the $t$-dependence in order to simplify the notation. In these new variables, equations \eqref{eqn:2.1}, \eqref{eqn:4.1} and \eqref{eqn:4.3} become
\begin{align*}
&\dot{x}=2x\eta,\\
&\frac{d^3}{d t^3}u+\frac{k}{t}\ddot{u}+\tr{(\eta)}\ddot{u}+\tr(\dot{\eta})\dot{u}-\frac{k}{t^2}\dot{u}-2\ddot{u}\dot{u}-\epsilon\dot{u}=0,\\
&\dot{\eta}=-\frac{k}{t^2}\MI_++\frac{1}{t^2}\MI_+-\frac{k}{t}\eta-\frac{1}{t}\tr{(\eta)}\MI_++\frac{1}{t}\dot{u}\MI_++r-\tr{(\eta)}\eta+\dot{u}\eta+\frac{\epsilon}{2}\MI.
\end{align*}
It is convenient to change again variables using
\begin{equation*}
y=x\eta,
\end{equation*}
so that we do not have to deal with the quadratic term $x\eta$. We then obtain
\begin{align}
\label{eqn:4.4}
&\dot{x}=2y,\\ \notag
&\frac{d^3}{d t^3}u+\frac{k}{t}\ddot{u}+\tr{(x^{-1}y)}\ddot{u}-\frac{k}{t^2}\dot{u}-2\tr{(x^{-1}yx^{-1}y)}\dot{u}+\tr{(x^{-1}\dot{y})}\dot{u}\\
\label{eqn:4.5}
&-2\ddot{u}\dot{u}-\epsilon\dot{u}=0,\\
\notag
&\dot{y}=(1-k)\frac{1}{t^2}x_+-\frac{k}{t}y-\frac{1}{t}\tr{(x^{-1}y)}x_++\frac{1}{t}\dot{u}x_++2yx^{-1}y+xr\\ 
\label{eqn:4.6}
&-\tr{(x^{-1}y)}y+\dot{u}y+\frac{\epsilon}{2}x,
\end{align}
with initial conditions on $y$ given by
\begin{equation*}
y_+(0)=0\quad\text{and}\quad y_-(0)=L_1(v_0).
\end{equation*}
At this point, we need the formula for the Ricci tensor of the homogeneous metric $g$ on the homogeneous space $P=G/K$, where $G$ is a compact Lie group. It has the following expression (see \cite{Besse:1} p. 185 for a derivation of this formula):
\begin{equation*}
\begin{aligned}
\Ric(X,Y)=&-\frac{1}{2}\tr_{\mathfrak g}(\ad(X)\ad(Y))-\frac{1}{2}\sum_{ij}{g([X,X_i]_{\mathfrak p},[Y,X_j]_{\mathfrak p})g^{ij}}\\
&+\frac{1}{4}\sum_{ijpq}{g(X,[X_i,X_p])g(Y,[X_j,X_q])g^{ij}g^{pq}},
\end{aligned}
\end{equation*}
for any basis $\{X_i\}_{i=1}^n$ of $\mathfrak p$ and for all $X,Y\in\mathfrak p$. Note that our expression for the Ricci tensor is simpler than $(7.38)$ in \cite{Besse:1}. This is due to the fact that $G$ is compact and hence unimodular.

The metric $\widehat{g}$ induces a $G$-invariant background metric $\widehat{g}_0$ on $NQ$. In fact, $\widehat{g}$ induces inner products on $\mathfrak p_-$, which can be identified with $T_qQ$, and on $V$, which can be identified with $N_qQ$.
Considering bases $\{U_\alpha\}_{\alpha=1}^k$ of $\mathfrak p_+$ and $\{Z_i\}_{i=k+1}^n$ of $\mathfrak p_-$, which are orthonormal with respect to the background metric $\widehat{g}_0$, the inverse of $g$ splits as follows
\begin{equation*}
g^{\alpha\beta}=\frac{1}{t^2}x_+^{\alpha\beta},\quad g^{ij}=x_-^{ij}.
\end{equation*}
Consequently, the Ricci endomorphism splits into a regular part and a singular part:
\begin{equation*}
r=\frac{1}{t^2}r_{\text{sing}}+r_{\text{reg}},
\end{equation*}
which are given in Lemma 3.1 of \cite{Wang-Eschenburg:1}.

We also have that
\begin{equation*}
\begin{aligned}
x_+r_+&=\frac{1}{t^2}g_+r_+=\frac{1}{t^2}\Ric_+,\\
x_-r_-&=g_-r_-=\Ric_-.
\end{aligned}
\end{equation*}
Hence, equation \eqref{eqn:4.6} becomes
\begin{equation*}
\dot{y}=\frac{1}{t^2}A(x)+\frac{1}{t}B(x,y)+C(x,y,t),
\end{equation*}
where
\begin{equation}\label{eqn:4.7}
\begin{aligned}
&A(x)=(1-k)x_++xr_{\text{sing}},\\
&B(x,y)=-ky-\tr{(x^{-1}y)}x_++\dot{u}x_+,\\
&C(x,y,t)=2yx^{-1}y+xr_{\text{reg}}-\tr{(x^{-1}y)}y+\dot{u}y+\frac{\epsilon}{2}x.
\end{aligned}
\end{equation}
We can now substitute $\dot{y}$ in \eqref{eqn:4.5} with the expression given by \eqref{eqn:4.6}. Then, Equation \eqref{eqn:4.5} becomes
\begin{equation*}
\frac{d^3}{d t^3}u=\frac{1}{t^2}\widetilde{A}(\dot{u})+\frac{1}{t}\widetilde{B}(\dot{u},\ddot{u})+\widetilde{C}(\dot{u},\ddot{u},t),
\end{equation*}
where
\begin{align*}
&\widetilde{A}(\dot{u})=k\dot{u}+(k-1)\tr{(x^{-1}x_+)}\dot{u}-\tr{(r_{\text{sing}})}\dot{u},\\
&\widetilde{B}(\dot{u},\ddot{u})=-k\ddot{u}+k\tr{(x^{-1}y)}\dot{u}+\tr{(x^{-1}y)}\tr{(x^{-1}x_+)}\dot{u}-\tr{(x^{-1}x_+)}\dot{u}^2\\
\begin{split}
&\widetilde{C}(\dot{u},\ddot{u},t)=-\tr{(x^{-1}y)}\ddot{u}-\tr{(r_{\text{reg}})}\dot{u}+\tr{(x^{-1}y)}\tr{(x^{-1}y)}\dot{u}-\tr{(x^{-1}y)}\dot{u}^2\\
&\quad\quad\quad\quad\quad\,-(n+1)\frac{\epsilon}{2}\dot{u}+2\ddot{u}\dot{u}+\epsilon\dot{u},
\end{split}
\end{align*}
are analytic functions.

We then obtain that the system \eqref{eqn:4.4}-\eqref{eqn:4.6}, with the above initial conditions, becomes the initial value problem given by
\begin{align*}
&\dot{x}=2y,\\
&\frac{d^3}{d t^3}u=\frac{1}{t^2}\widetilde{A}(\dot{u})+\frac{1}{t}\widetilde{B}(\dot{u},\ddot{u})+\widetilde{C}(\dot{u},\ddot{u},t),\\
&\dot{y}=\frac{1}{t^2}A(x)+\frac{1}{t}B(x,y)+C(x,y,t),\\
&x(0)=\MI,\\
&y_+(0)=0\quad\text{and}\quad y_-(0)=L_1(v_0),\\
&\dot{u}(0)=0.
\end{align*}
Therefore, the initial value problem for cohomogeneity one Ricci solitons has been reduced to an initial value problem for a system of non linear ordinary differential equations of order one in $x,y$ and of order two in $\dot{u}$ with a singular point at the origin.

By Lemma 2.2 of \cite{Dancer-Wang:2}, Ricci solitons are real analytic. Therefore, we can solve the system by applying the method of asymptotic power series, which is described in \cite{Wasow:1}, Chapter 9. This method consists, first of all, of showing that there always exists a formal power series solution of an appropriate type. Then, after having a formal power series solution, one can apply Theorem 7.1 of \cite{Malgrange:1} to get a genuine solution.

We will see that the reason why there always exists a formal power series solution, is due to the geometric nature of the equations.

\section{Solution to the initial value problem}
The initial value problem considered in Section \ref{sec:4} has the following general form
\begin{align}
&\dot{x}=2y,\label{eqn:5.1}\\
&\frac{d^3}{d t^3}u=\frac{1}{t^2}\widetilde{A}(\dot{u})+\frac{1}{t}\widetilde{B}(\dot{u},\ddot{u})+\widetilde{C}(\dot{u},\ddot{u},t)\label{eqn:5.2}\\
&\dot{y}=\frac{1}{t^2}A(x)+\frac{1}{t}B(x,y)+C(x,y,t),\label{eqn:5.3}\\
&x(0)=a,\label{eqn:5.4}\\
&y(0)=b,\label{eqn:5.5}\\
&\dot{u}(0)=0,\label{eqn:5.6}
\end{align}
where $x(t),y(t),a,b\in \Sym^2(V\oplus\mathfrak p_-)^K$, $u(t)$ is smooth function and $A$, $B$, $C$, $\widetilde{A}$, $\widetilde{B}$ and $\widetilde{C}$ are analytic functions.

As the left-hand sides of \eqref{eqn:5.2} and \eqref{eqn:5.3} do not have $\frac{1}{t}$ or $\frac{1}{t^2}$ terms, $A$ and $\widetilde{A}$ must satisfy the following initial conditions
\begin{align}
&A(x(0))=A(a)=0\quad\text{and}\quad2(d A)_a\cdot b+B(a,b)=0,\label{eqn:5.12}\\
&\widetilde{A}(\dot{u}(0))=\widetilde{A}(0)=0\quad\text{and}\quad(d\widetilde{A})\vert_{t=0}\ddot{u}(0)+\widetilde{B}(0,\ddot{u}(0))=0.\label{eqn:5.13}
\end{align}
We want to show that there always exists a formal power series solution. So let
\begin{equation*}
x(t)=\sum_{m=0}^\infty{\frac{x_m}{m!}t^m},\quad y(t)=\sum_{m=0}^\infty{\frac{y_m}{m!}t^m},
\end{equation*}
with
\begin{equation*}
x_{m+1}=2y_m,\quad\forall m\geq0,
\end{equation*}
and
\begin{equation*}
u(t)=\sum_{m=0}^\infty{\frac{u_m}{m!}t^m}.
\end{equation*}
Then, let
\begin{equation*}
A(x(t))=\sum_{m=0}^\infty{\frac{A_m}{m!}t^m},\quad B(x(t),y(t))=\sum_{m=0}^\infty{\frac{B_m}{m!}t^m},\quad C(x(t),y(t),t)=\sum_{m=0}^\infty{\frac{C_m}{m!}t^m},
\end{equation*}
and
\begin{equation*}
\widetilde{A}(\dot{u}(t))=\sum_{m=0}^\infty{\frac{\widetilde{A}_m}{m!}t^m},\quad\widetilde{B}(\dot{u}(t),\ddot{u}(t))=\sum_{m=0}^\infty{\frac{\widetilde{B}_m}{m!}t^m},\quad\widetilde{C}(\dot{u}(t),\ddot{u}(t),t)=\sum_{m=0}^\infty{\frac{\widetilde{C}_m}{m!}t^m}.
\end{equation*}
Substituting the above expressions in \eqref{eqn:5.2} and \eqref{eqn:5.3}, respectively, we get
\begin{align}
\frac{1}{2}x_{m+2}&=\frac{A_{m+2}}{(m+2)(m+1)}+\frac{B_{m+1}}{m+1}+C_m,\label{eqn:5.8}\\
u_{m+3}&=\frac{\widetilde{A}_{m+2}}{(m+2)(m+1)}+\frac{\widetilde{B}_{m+1}}{m+1}+\widetilde{C}_m.\label{eqn:5.9}
\end{align}
By definition, we have that
\begin{align*}
A_{m+2}&=\frac{d^{m+2}}{d t^{m+2}}(A(x(t))\Big\vert_{t=0}=\frac{d^{m+1}}{d t^{m+1}}\left(\frac{d}{d t}(A(x(t)))\right)\Big\vert_{t=0}\\
&=\frac{d^{m+1}}{d t^{m+1}}(d_{x(t)}A\cdot\dot{x}(t))\Big\vert_{t=0}\\
&\equiv(d_{x(t)} A)_a\cdot x_{m+2}\quad(\text{mod}\,x_1,\dots,x_{m+1}),\\
B_{m+1}&=\frac{d^{m+1}}{d t^{m+1}}(B(x(t),y(t)))\Big\vert_{t=0}=\frac{d^m}{d t^m}\left(\frac{d}{d t}(B(x(t),y(t)))\right)\Big\vert_{t=0}\\
&=\frac{d^m}{d t^m}(\partial_{x(t)}B\cdot\dot{x}(t)+\partial_{y(t)}B\cdot\dot{y}(t)+x_+\ddot{u})\Big\vert_{t=0}\\
&\equiv\frac{1}{2}(\partial_{y(t)}B)_{(a,b)}\cdot x_{m+2}\quad(\text{mod}\,x_1,\dots,x_{m+1}),\\
C_m&\equiv0\quad(\text{mod}\,x_1,\dots,x_{m+1}).
\end{align*}
Using the same strategy, we also have that
\begin{align*}
\widetilde{A}_{m+2}&\equiv(d_{\dot{u}(t)}\widetilde{A})\vert_{t=0}u_{m+3}\quad(\text{mod}\,u_1,\dots,u_{m+2}),\\
\widetilde{B}_{m+1}&\equiv(\partial_{\ddot{u}(t)}\widetilde{B})\vert_{t=0}u_{m+3}\quad(\text{mod}\,u_1,\dots,u_{m+2}),\\
\widetilde{C}_m&\equiv0\quad(\text{mod}\,u_1,\dots,u_{m+2}).
\end{align*}
Hence, equations \eqref{eqn:5.8} and \eqref{eqn:5.9} become
\begin{align}
&x_{m+2}=2\frac{(d_{x(t)} A)_a\cdot x_{m+2}}{(m+2)(m+1)}+\frac{(\partial_{y(t)}B)_{(a,b)}\cdot x_{m+2}}{m+1}+\frac{D_m}{m+1},\label{eqn:5.14}\\
&u_{m+3}=\frac{(d_{\dot{u}(t)}\widetilde{A})\vert_{t=0}u_{m+3}}{(m+2)(m+1)}+\frac{(\partial_{\ddot{u}(t)}\widetilde{B})\vert_{t=0} u_{m+3}}{m+1}+\frac{\widetilde{D}_m}{m+1},\label{eqn:5.15}
\end{align}
for some functions $D_m$ of $x_1,\dots,x_{m+1}$ and $\widetilde{D}_m$ of $u_1,\dots,u_{m+2}$. If we now define the following two operators
\begin{align}
&\mathcal L_{m}=(m+1)\MI-\frac{2}{m+2}(d_{x(t)} A)_a-(\partial_{y(t)}B)_{(a,b)},\label{eqn:5.11}\\
&\mathcal{\widetilde{L}}_m=(m+1)-\frac{1}{m+2}(d_{\dot{u}(t)}\widetilde{A})\vert_{t=0}-(\partial_{\ddot{u}(t)}\widetilde{B})\vert_{t=0},\notag
\end{align}
we need to have that
\begin{equation*}
\mathcal L_{m}\cdot x_{m+2}=D_m\quad\text{and}\quad\mathcal{\widetilde{L}}_mu_{m+3}=\widetilde{D}_m,
\end{equation*}
which give necessary and sufficient conditions to the existence of a formal power series solution:
\begin{equation}\label{eqn:5.10}
D_m\in\image(\mathcal L_{m})\quad\text{and}\quad\widetilde{D}_m\in\image(\mathcal{\widetilde{L}}_m),
\end{equation}
for all $m\geq0$. We have that $\mathcal L_m$ and $\mathcal{\widetilde{L}}_m$ are invertible for all $m\geq m_0$, for some $m_0$. In fact, $\mathcal{\widetilde{L}}_m$ is bounded and $d A$ and $d B$ are bounded as well. For this reason, if $m$ is large, $\mathcal L_m$ is close to a multiple of the identity and hence invertible. This implies that, if \eqref{eqn:5.10} is satisfied for $m<m_0$, we can fix further initial conditions, namely $x_1,\dots,x_{m_0}$ and $u_1,\dots,u_{m_0}$ satisfying equations \eqref{eqn:5.14} and \eqref{eqn:5.15} respectively, such that the formal power series solution is uniquely determined. 

As we said before, $a$ and $b$ are two $K$-invariant endomorphisms which preserve the splitting of $\mathfrak p$ in a $+$ part and in a $-$ part. Moreover, as we saw in section \ref{sec:4}, they are given by
\begin{equation*}
a=\MI,\quad b_+=0\quad\text{and}\quad b_-=L_1(v_0),
\end{equation*}
where we have that $L_1\in W_1^-$, which implies that $\tr{(L_1)}=0$.

Now, using expressions given in \eqref{eqn:4.7} and Lemmas 4.2 and 4.4 of \cite{Wang-Eschenburg:1}, we can write the operator $\mathcal L_m$ as follows. First of all, we have that
\begin{align*}
&(d_{x(t)} A)_a\cdot\xi=(d_{x(t)} r_{\text{sing}})_{\MI}\cdot\xi,\\
&B(a,b)=-kb,\\
&(\partial_yB)_{(a,b)}\cdot\xi=-k\xi-\tr{(\xi)}\MI_+,\\
&C(x,y,t)=2yx^{-1}y+xr_{\text{reg}}-\tr{(x^{-1}y)}y+\dot{u}y+\frac{\epsilon}{2}x,
\end{align*}
where $\xi\in \Sym^2(V\oplus\mathfrak p_-)^K$ and 
\begin{align*}
(d_{x(t)} r_{\text{sing}})_{\MI}\cdot\xi_+&=(k+1)\xi_+-2\tr{(\xi_+)}\MI_+,\\
(d_{x(t)} r_{\text{sing}})_{\MI}\cdot\xi_-&=\frac{1}{2}\mathcal C\cdot\xi_-,
\end{align*}
where $\mathcal C$ is an operator defined by $\mathcal C=-\sum_{\alpha=1}^k{\ad(U_\alpha)^2}$. Note that, if the $H$-homogeneous standard metric on $S^k$ is normal, we have a bi-invariant metric on $\mathfrak h$ and we can extend the basis $\{U_\alpha\}$ of $\mathfrak p_+\subset\mathfrak h$ to an orthonormal basis $\{V_\alpha\}$ of $\mathfrak h$, equipped with this bi-invariant metric, and we have that $\mathcal C=-\sum_\alpha{\ad(V_\alpha)^2}$ is the Casimir operator for the adjoint representation on $\mathfrak p_-$ and $\End(\mathfrak p_-)$. Note that we obtained, apart from $C$, the same expressions as \cite{Wang-Eschenburg:1} (p. 129).

Substituting these expression in \eqref{eqn:5.11}, we have that
\begin{equation*}
\mathcal L_m\cdot\xi=(m+1)\xi-\frac{2}{m+2}(d r_{\text{sing}})_{\MI}\cdot\xi+k\xi+\tr{(\xi)}\MI_+.
\end{equation*}
Furthermore, by Lemma 4.6 of \cite{Wang-Eschenburg:1}, we have that
\begin{align*}
&(\mathcal L_m\cdot\xi)_+=m\left(1+\frac{k+1}{m+2}\right)\xi_++\left(\frac{4}{m+2}\tr{(\xi_+)}+\tr{(\xi)}\right)\MI_+,\\
&(\mathcal L_m\cdot\xi)_-=(m+1+k)\xi_--\frac{1}{m+2}\mathcal C\cdot\xi_-.
\end{align*}
Again by Lemmas 4.2 and 4.4 of \cite{Wang-Eschenburg:1} and by the above expressions, we also have that
\begin{align*}
(d_{\dot{u}(t)}\widetilde{A})\vert_{t=0}f&=kf+(k-1)\tr{(x^{-1}(0)x_+(0))}f-\tr{(r_{\text{sing}}(0))}f\\
&=kf+(k-1)kf-\tr{((k-1)\MI_+)}f\\
&=kf+(k-1)kf-(k-1)kf\\
&=kf,\\
(\partial_{\ddot{u}(t)}\widetilde{B})\vert_{t=0}f&=-kf,
\end{align*}
so that the operator $\mathcal{\widetilde{L}}_m$ becomes
\begin{equation}\label{eqn:5.16}
\mathcal{\widetilde{L}}_m=(m+1)-\frac{k}{m+2}+k.
\end{equation}
We now have to verify that the initial conditions \eqref{eqn:5.12} and \eqref{eqn:5.13} hold and that \eqref{eqn:5.10} is satisfied. Note that \eqref{eqn:5.12} is satisfied as explained in \cite{Wang-Eschenburg:1} (p. 130). Then, we have that $\widetilde{A}(\dot{u}(0))=0$, because $\dot{u}(0)=0$. Moreover, 
\begin{equation*}
\begin{aligned}
(d_{\dot{u}(t)}\widetilde{A})\vert_{t=0}\ddot{u}(0)=&k\ddot{u}(0)+(k-1)\tr(x^{-1}(0)x_+(0))\ddot{u}(0)-\tr{(r_{\text{sing}}(0))}\ddot{u}(0)\\
=&k\ddot{u}(0)
\end{aligned}
\end{equation*}
by Lemmas 4.2 and 4.4 of  \cite{Wang-Eschenburg:1} and $\widetilde{B}(0,\ddot{u}(0))=-k\ddot{u}(0)$ by definition. Hence, \eqref{eqn:5.13} holds as well. As the initial conditions hold, we need to verify equation \eqref{eqn:5.10} and to show that 
\begin{equation*}
x^l(t):=\sum_{m=0}^l{\frac{x_m}{m!}t^m}
\end{equation*}
defines a smooth $G$-invariant metric around $Q$ for all $l$. By Lemma \ref{Lemma:3.1} this means that we need to prove that $x_m\in\text{ev}(W_m)$ for all $m$.

In  \cite{Wang-Eschenburg:1}, the authors show that $\mathcal L_m(\text{ev}(W_{m+2}))\subset\text{ev}(W_m)$, by decomposing $\text{ev}(W_{m+2})$ into eigenspaces of $\mathcal L_m$ and showing that the only eigenspaces corresponding to nonzero eigenvalues lie in $\text{ev}(W_m)\subset\text{ev}(W_{m+2})$, (we can always modify the degree of a homogeneous polynomial by an even factor without changing its value on the sphere). Moreover, $\mathcal L_m$ maps $\text{ev}(W_m)$ bijectively onto itself if $m>0$. Hence, the equation $\mathcal L_m\cdot x_{m+2}=D_m$ has a solution if and only if $D_m\in\text{ev}(W_m)$. Moreover, we also have that the kernel of $\mathcal L_m$ is isomorphic to $W^-_{m+2}/W^-_m$. We need now to show that $x_p\in\text{ev}(W_p)$ for all $p$ and that $D_m\in\text{ev}(W_m)$ for all $m$. In  \cite{Wang-Eschenburg:1}, the authors show this by induction over $m$. They show that there exists a solution $x_{m+2}$ of $\mathcal L_m\cdot x_{m+2}=D_m$, but we can add an arbitrary element of the kernel of the operator considered, as it is not trivial. 

Considering the multiplicative operator defined by \eqref{eqn:5.16}, we can see that it maps $\text{ev}(W_{m+3})$ to itself. So we need to show that $\widetilde{D}_m\in\text{ev}(W_{m+3})$ and that 
\begin{equation*}
u^m(t):=\sum_{p=0}^{m}{\frac{u_p}{p!}t^p}
\end{equation*}
defines a smooth $G$-invariant function around $Q$ for all $m$. This means that we have to show that $u_p\in\text{ev}(W_p)$ for all $p$. We can prove that $\widetilde{D}_m\in\text{ev}(W_{m+3})$ by induction over $m$. We have that $u_1=0\in\text{ev}(W_1)$. Then, suppose that $u_p\in\text{ev}(W_p)$ for $p=2,\dots,m+2$, which implies that $u^{m+2}(t)$ is even in $t$, and consider
\begin{equation*}
\widehat{u}(t)=u^{m+2}(t)=\sum_{p=0}^{m+2}{\frac{u_p}{p!}t^p},
\end{equation*}
which defines a smooth $G$-invariant function around $Q$, because of the discussion in Section \ref{section 3}. By definition, we have that $\widehat{u}$ satisfies equation \eqref{eqn:5.2}. Let $\widehat{A},\widehat{B}$ and $\widehat{C}$ be the the analogues of $\widetilde{A},\widetilde{B},\widetilde{C}$ for $\widehat{u}$. Moreover, let $\widehat{D}_m$ be some function of $\widehat{u}_1,\dots,\widehat{u}_{m+2}$ which satisfies an  analogue of equation \eqref{eqn:5.15} for $\widehat{u}$. Now, as $\widehat{u}_{m+3}=0\in\text{ev}(W_{m+3})$ and $\mathcal{\widetilde{L}}_m\widehat{u}_{m+3}=\widehat{D}_m$, we have that $\widehat{D}_m=0$. Then, by equations \eqref{eqn:5.9} and \eqref{eqn:5.15} and by recalling the expression of $\widetilde{C}$, we have that
\begin{equation*}
\frac{\widetilde{D}_m}{m+1}=\frac{\widetilde{D}_m-\widehat{D}_m}{m+1}=\widetilde{C}_m-\widehat{C}_m=0\in\text{ev}(W_{m+3}).
\end{equation*}
Hence, $\widetilde{D}_m\in\text{ev}(W_{m+3})$ and a solution $u_{m+3}$ to $\mathcal{\widetilde{L}}_mu_{m+3}=\widetilde{D}_m$ exists in $\text{ev}(W_{m+3})$.

The indeterminacy is just in the operator $\mathcal L_m$ and it is the same as in the Einstein case. In \cite{Wang-Eschenburg:1}, the authors describe this indeterminacy in the formal power series solution. If $m>0$, after solving the $-$ part of $\mathcal L_m\cdot x_{m+2}=D_m$, the $+$ part is uniquely determined. On the contrary, if $m=0$, the trace free part of $(x_2)_+$ is arbitrary. This is to be expected, as the values of $x(0)$ and $y_+(0)$ are fixed by the geometry of the problem and this implies that the usual freedom in the initial value problem lies in the trace free part of $(x_2)_+$. Furthermore, by section 1 of \cite{Wang-Eschenburg:1}, we see that the spaces $W_m^-$ eventually stabilise. So, suppose that
\begin{equation*}
W^-_{2m}=W^-_{2m_0}\quad\text{and}\quad W^-_{2m+1}=W^-_{2m_1+1},
\end{equation*}
for all $m>m_0$ and for all $m>m_1$, respectively. Hence, as $\ker(\mathcal L_m)\backsimeq W^-_{m+2}/W^-_m$, we have that the indeterminacy of the initial value problem considered is given by
\begin{equation*}
(W^-_{2m_0}/W^-_0)\oplus(W^-_{2m_1+1}/W^-_1),
\end{equation*}
where, in particular,
\begin{equation*}
W_0=\Sym^2(V)^H\oplus \Sym^2(\mathfrak p_-)^H.
\end{equation*}
The formal solution to the initial value problem \eqref{eqn:5.1}-\eqref{eqn:5.6} has the property that if truncated at any order it gives a smooth $G$-invariant metric and a smooth $G$-invariant function on $E$. Now, by Theorem 7.1 in \cite{Malgrange:1}, we obtain a genuine solution to the problem considered which is defined on a small interval $[0,T]$. The genuine solution may also be obtained by carrying out the Picard iteration directly, as shown in \cite{Wang-Eschenburg:1}. From this, one can see that $x_m(t)$ defines a smooth $G$-invariant metric on a tubular neighbourhood of radius $T$ around $Q$. Moreover, by Lemma 1.2 of \cite{Wang-Eschenburg:1}, we can choose $m$ to be at least 3. Hence, the solution $x$ gives a $C^3$ $G$-invariant metric which satisfies the appropriate equation. Finally, by Lemma 2.2 of \cite{Dancer-Wang:2}, the pair $(x,u)$ gives a smooth $G$-invariant Ricci soliton on the tubular neighbourhood of radius $T$ around the singular orbit $Q$.

To conclude this section, as the reader may not be familiar with \cite{Wang-Eschenburg:1}, we describe some examples, which show how to compute the indeterminacy explicitly.
\begin{example}[\cite{Wang-Eschenburg:1}]
Let
\begin{equation*}
G=SO(p+n),\quad H=SO(p)\times SO(n),\quad\text{and}\quad K=SO(p)\times SO(n-1).
\end{equation*}
Then, the manifold we are dealing with has dimension $np+n+1$. $H$ acts effectively and transitively on the unit sphere in $V$. By Lemma 1.2 of \cite{Wang-Eschenburg:1}, we have that $W_m^+\backsimeq\Hom(\Sym^m(V),\Sym^2(V))^H$ is isomorphic to zero if $m$ is odd and that all these spaces are isomorphic if $m$ is even. Hence, we can compute the indeterminacy in $(x_2)_+$, which is given by the dimension of $\Hom(\Sym^2(V),\Sym^2(V))^H$, which is 1, by Lemma 1.2 of \cite{Wang-Eschenburg:1}. Now, we need to compute the dimension of $W_m^-\backsimeq\Hom(\Sym^m(V),\Sym^2(\mathfrak p_-))^H$. As an $H$-representation $\mathfrak p_-$ is the tensor product of the standard representation $\rho_p$ of $SO(p)$ and the standard representation $\rho_n$ of $SO(n)$. On the other hand, $V$ is the tensor product of the trivial representation $\boldsymbol 1$ of $SO(p)$ and $\rho_n$. It is well known (see \cite{Fulton-Harris:1} p. 296) that 
\begin{equation*}
\Sym^m(V)=\sigma_m\oplus\sigma_{m-2}\oplus\cdots\oplus\sigma_{m-2[\frac{m}{2}]},
\end{equation*}
where $\sigma_k$ is the irreducible representation of $SO(n)$ with dominant weight $k$ times that of $\rho_n$. We also have that
\begin{equation*}
\Sym^2(\mathfrak p_-)=(\Sym^2(\rho_p)\otimes \Sym^2(\rho_n))\oplus(\Lambda^2(\rho_p)\otimes\Lambda^2(\rho_n)),
\end{equation*}
which decomposes as
\begin{equation*}
(\sigma_2\otimes\sigma_2)\oplus(\boldsymbol 1\otimes\sigma_2)\oplus(\sigma_2\otimes\boldsymbol 1)\oplus(\boldsymbol 1\otimes\boldsymbol 1)\oplus(\ad_p\otimes\ad_n),
\end{equation*}
where we used the fact that the adjoint representation of $SO(n)$ is equivalent to $\Lambda^2V$.

Applying Schur's Lemma, we have that the dimension of $W_m^-$ is 2 when $m\geq2$ is even, it is zero if $m$ is odd. Finally, $\dim(W_0^-)=\dim(\Sym^2(\mathfrak p_-)^H)=1$. Thus, the indeterminacy which occurs with $(x_2)_-$ has dimension 1. We have that the choice of the initial metric is unique up to homothety and the only choice for the initial shape operator is zero, because $V$ is not a summand in $\Sym^2(\mathfrak p_-)$. This means that the singular orbit must be totally geodesic. Anyway, there is a one-dimensional freedom in choosing $(x_2)_+$.
\end{example}
\begin{example}
Let 
\begin{equation*}
G=SO(n+2),\quad K=SO(n),\quad H=SO(n+1).
\end{equation*}
The manifold $M$ on which $G$ acts has dimension $2n+2$, the principal orbits for this action are given by the Stiefel manifold, which is the homogeneous manifold $SO(n+2)/SO(n)$ for $n\geq2$, and the singular orbit is given by the sphere $S^{n+1}$, which has codimension $n+1$ in $M$.

We want to compute the indeterminacy in the initial value problem for cohomogeneity one gradient Ricci solitons. We have that $V\backsimeq\mathbb R^{n+1}$, as $H$-representation, is given by the standard orthogonal representation $\rho_{n+1}$. Similarly, $\mathfrak p_-\backsimeq\mathbb R^{n+1}$, as $H$-representation, is given by $\rho_{n+1}$. Hence,
\begin{equation*}
\Sym^2(\mathfrak p_-)=\Sym^2(\rho_{n+1})=\sigma_2\oplus\boldsymbol 1.
\end{equation*}
Note that the assumption \eqref{assumption} is satisfied, because we assume that the metric $g_t$ is diagonal with respect to the following decomposition:
\begin{equation*}
\mathfrak g=\mathfrak k\oplus\underbrace{\mathfrak p_1\oplus\mathfrak p_2\oplus\mathfrak p_3}_{\mathfrak p},
\end{equation*}
where $\mathfrak p_1$ and $\mathfrak p_2$ are n-dimensional $SO(n)$-representations and $\mathfrak p_3$ is the trivial $SO(n)$-representation. With this decomposition we have that $\mathfrak p_-=\mathfrak p_2\oplus\mathfrak p_3$. By Schur's Lemma we then have that the dimension of $W_m^-$ is zero if $m$ is odd and it is one if $m\geq2$ is even. Moreover, $W_0^-$ has dimension one. So, we have that the indeterminacy in $(x_2)_-$ is zero. We can now compute the dimension of $W_m^+$, which is zero if $m$ is odd and 1 if either $m$ is zero or $m\geq2$ is even. So the indeterminacy in $(x_2)_+$, which is given by the dimension of $W_2^+$, is one. So we just have a one-dimensional freedom in choosing $(x_2)_+$.
\newline

Instead of considering $H=SO(n+1)$, we could consider $H=SO(2)\times SO(n)$, so that the singular orbit $Q=G/H$ has codimension 2. We have that $V\backsimeq\mathbb R^2$ and, as an $H$-representation, it is given by $\rho_2\otimes\boldsymbol 1$. On the other hand, $\mathfrak p_-\backsimeq\mathbb R^{2n}$, as an $H$-representation, is given by $\boldsymbol 1\otimes\rho_n\otimes\rho_n$. We then have that
\begin{equation*}
\begin{aligned}
\Sym^2(\mathfrak p_-)&=\Sym^2(\rho_n\otimes\rho_n)=(\Sym^2(\rho_n)\otimes\Sym^2(\rho_n))\oplus(\Lambda^2(\rho_n)\otimes\Lambda^2(\rho_n))\\
&=(\sigma_2\otimes\sigma_2)\oplus(\boldsymbol 1\otimes\sigma_2)\oplus(\sigma_2\otimes\boldsymbol 1)\oplus(\boldsymbol 1\otimes\boldsymbol 1)\oplus(\ad_n\otimes\ad_n).
\end{aligned}
\end{equation*}
Then we have that the dimension of $W_m^-$ is two when $m\geq2$ is even and it is zero when $m$ is odd. Moreover, the dimension of $W_0^-$ is one. So we have that the indeterminacy in $(x_2)_-$ is one. As the dimension of $W_2^+$ is one, the indeterminacy in $(x_2)_+$ is one.
\end{example}

\section*{Acknowledgements}
I would like to thank my supervisor Prof. Andrew Dancer for his useful comments and his supervision in writing the paper and Thomas Madsen for interesting discussions.

\bibliographystyle{amsplain}
\bibliography{citazioni}

\end{document}